\documentclass{article}
\usepackage{amssymb,amsmath}

\newtheorem{theorem}{Theorem}[section]

\newtheorem{lemma}[theorem]{Lemma}

\newcommand{\zz}{{\mathbb{Z}}}
\newcommand{\ff}{{\mathbb{F}}}

\title{ Orthogonal and Symplectic Black Box Groups,\\ Revisited}

\author{Alexandre V. Borovik\thanks{Research supported by the Royal Society and Leverhulme Trust.}}

\begin{document}

\maketitle

\pagestyle{myheadings}
\markright{{\scriptsize A. V. Borovik $\bullet$ Orthogonal and Symplectic Groups, Revisited $\bullet$ 21.10.01}}

\begin{abstract}
We propose a simple one sided Monte-Carlo algorithm to
distinguish, to any given degree of certainty, the symplectic
group $C_n(q) = {\rm PSp}_{2n}(q)$ from the orthogonal group
$B_n(q) = \Omega_{2n+1}(q)$ where $q>3$ is odd  and $n$ and $q$
are given. The algorithm does not use an order oracle and works in
polynomial, of $n\log q$, time.

\medskip
This paper corrects an error in the previously published version
of the algorithm \cite{old}.

\medskip

\end{abstract}

\section{Introduction}

The aim of this paper is to correct an error in the previously
published version \cite{old} of a  Monte-Carlo algorithm which
distinguishes, to any given degree of certainty, the symplectic
group $C_n(q) = {\rm PSp}_{2n}(q)$ from the orthogonal group
$B_n(q) = \Omega_{2n+1}(q)$ where $q>3$ is odd and $n$ and $q$ are
given. Namely we prove the following theorem.

\begin{theorem}
\label{main} There is a polynomial time one sided Monte-Carlo
algorithm which, when given a black box group $G$ isomorphic to a
simple classical group of type $B_n$ or\/ $C_n$ defined over some
field of given size $q>3$ and odd characteristic, can decide
whether $G$ is isomorphic to $C_n$ or not. If $G$ is  $C_n$ then,
with probability $ > 1/2$, the algorithm produces a witness
against isomorphism $G \simeq B_n(q)$.
\end{theorem}

The authors would like to thank Eamonn O'Brien and Gunter Malle who found an error
in the original version of the algorithm (Lemma~2.4 of \cite{old}
was false). I am also grateful to Sukru Yalcinkaya who corrected
a considerable number of typographic errors in the text.

We use the notation and results of \cite{old}.

\subsection{Description of the algorithm: the easy case, $n$ is odd.}
Recall that $\Omega_{2n+1}(q)$ denotes the kernel of the spinor norm
$${\rm SO}_{2n+1}(q) \longrightarrow \{\,\pm 1\,\};
$$
it is a subgroup of index $2$ in ${\rm SO}_{2n+1}(q)$ and coincides with the commutator of
${\rm SO}_{2n+1}(q)$ if $n \geqslant 2$ and $q \geqslant 3$ (for details see, for example,
Artin \cite{artin}).
It is well known that ${\rm PSp}_2(q) \simeq \Omega_3(q) \simeq {\rm PSL}_2(q)$ and
${\rm PSp}_4(q) \simeq {\Omega}_5(q)$. Hence it is obviously enough to
give an algorithm for $n \geqslant 3$.

The easy case of the algorithm deals with the situation when $n$ is odd, and the hard
case is concerned with even $n$.

If $x$ is an element of even order, we use notation ${\rm i}(x)$ for the involution
in the cyclic group $\langle x\rangle$.

The structure of the algorithm in the easy case is the same as
described in \cite{old} and very simple.

\begin{description}

\item[{\rm INPUT:}] A black box group $G$ isomorphic to ${\rm PSp}_{2n}(q)$ or
${\Omega}_{2n+1}(q)$, $n\geqslant 3$ and $q>3$ given, $n$ is odd.

\item[{\rm OUTPUT:}] If $G$ is   ${\rm PSp}_{2n}(q)$ then, with probability
$ > 1/2$, the algorithm produces a witness against isomorphism  $G \simeq {\Omega}_{2n+1}(q)$.

\end{description}

\begin{itemize}
\item[(1)] Find an element $x$ of {\it good order} (the definition is given
below).  As it was discussed in \cite{old}, this can be done in
polynomial time.

\item[(2)] Compute the involution $i={\rm i}(x)$ in the cyclic group generated by
  $x$, this also can be done in polynomial time.

\item[(3)] {\em This step is vacuous in the easy case.}

\item[(4)] Compute $y:=i^gi^h$ for random elements $g,h \in G$. If
  $o(y)$ divides $q(q+1)$ or $q(q-1)$ then go to Step 5. Otherwise return
  {\it `$G$ is  ${\rm PSp}_{2n}(q)$' }, with the pair $(i^g, i^h)$ being a
  witness against the isomorphism
$G \simeq {\Omega}_{2n+1}(q)$.

\item[(5)] Repeat Step 4 at least  $20n$ times. If $o(y)$ always divides $q(q+1)$ or
  $q(q-1)$ then return {\it `$G$ is not\/ ${\rm PSp}_{2n}(q)$'}. The error probability is $< 1/2$
  in this case.
\end{itemize}

The number of necessary repetitions in  Step 5 is a very crude
estimate, which was sufficient for our theoretical purposes. In
practice, if $G$ is of type $C_n$ and $q>3$, then one does get
an  element $i^g\cdot i^h$ whose order does not divide $q(q+1)$ or
$q(q-1)$ at the very first iteration.

\subsection{Description of the algorithm: the hard case, $n$ is even.}

Let $p$ be a prime number  $p$, $q = p^k$  and $m \geqslant 2$.  A
prime dividing $q^m-1$ but not $q^i-1$ for $1 \leqslant i < m$ is
called a {\it primitive prime divisor of $q^m-1$}.  It will be
convenient to say that $p$ is a primitive prime divisor of
$q^m+1$ if $p$ is a primitive prime divisor of $q^{2m}-1$, see
Lemma\ref{lm:ppd(q)} for further discussion of primitive prime
divisors. In our case primitive prime divisors always exist by
Zsigmondy \cite{zs}. The {primitive divisor rank\/} ${\rm pdrank}
(g)$ is the maximal integer $k$ such that $o(g)$ is divisible by a
primitive prime divisor of $q^k-1$.

\begin{description}

\item[{\rm INPUT:}] A black box group $G$ isomorphic to
${\rm PSp}_{2n}(q)$ or ${\Omega}_{2n+1}(q)$, $n \geqslant 3$ and $q$ given, $n$ is even
(and hence $q^n \equiv 1 \bmod 4$).

\item[{\rm OUTPUT:}] If $G$ is  ${\rm PSp}_{2n}(q)$ then, with probability
$ > 1/2$, the algorithm constructs a witness against
the isomorphism  $G \simeq {\Omega}_{2n+1}(q)$.

\end{description}

\begin{itemize}
\item[(1)] Find an element $x$ of {\it good order} as in the easy case.

\item[(2)] Compute the involution $i={\rm i}(x)$ in the cyclic group generated by
  $x$.

\medskip

\begin{table}[p]
\begin{center}
\begin{tabular}{|| c |c | c ||}\hline && \\[1ex]

Involution $t_i$ $\backslash$ Group $G$ & ${\rm PSp}_{2n}(q)$ & $\Omega_{2n+1}(q)$\\[2ex] \hline
&& \\[1ex]

$t_1$ & ${\rm Sp}_{2}(q)\circ {\rm Sp}_{2n-2}(q)$
    & $\mathbb{Z}_{q\pm 1} \times \Omega_{2n-1}(q)$ \\[2ex] \hline
&& \\[1ex]
$t_i$ & ${\rm Sp}_{2i}(q)\circ {\rm Sp}_{2n-2i}(q)$ &
$\Omega^+_{2i}(q) \times \Omega_{2(n-i)+1}(q)$ \\[2ex]
$2 \leqslant i \leqslant n/2$ (in ${\rm PSp}_n(q)$)&  & ($i$ even)\\[2ex]
$2 \leqslant i < n-1$ (in $\Omega_{2n+1}(q)$)&& $\Omega^-_{2i}(q) \times \Omega_{2(n-i)+1}(q)$\\[2ex]
&& ($i$ odd  and $q^i \equiv -1 \bmod{4}$)\\[2ex]\hline
&& \\[2ex]
$t_{n-1}$ (in $\Omega_{2n+1}(q)$)&& $\Omega^+_{2n-2}(q) \times {\rm PSL}_3(q)$ \\[2ex]
&  & ($n$ is odd)\\[2ex]
&& \\[2ex]

&&$\Omega^-_{2n-2}(q) \times {\rm PSL}_3(q)$\\[2ex]
&& ($n$ is even  and $q^{n-1} \equiv -1 \bmod{4}$)\\[2ex]\hline

&&\\[2ex]
$t_n$ & $\frac{1}{2}{\rm GL}_n(q)$ & $\Omega^+_{2n}(q)$\\[2ex]
$q\equiv 1 \bmod{4}$ &&\\[2ex]\hline
&&\\[2ex]
$t_n$ & $\frac{1}{2}{\rm GU}_n(q^2)$ & $\Omega^-_{2n}(q)$\\[2ex]
$q\equiv -1 \bmod{4}$ &&\\[2ex]
&&\\[1ex]\hline
\end{tabular}
\end{center}
\caption{Involutions in groups ${\rm PSp}_{2n}(q)$ and
$\Omega_{2n}(q)$, $n$ even. The second and the third column show
the generalised Fitting subgroup $F^*(C_G(t_i))$  of $C_G(t_i)$
in the groups ${\rm PSp}_{2n}(q)$ and $\Omega_{2n+1}(q)$,
correspondingly. In all cases $C_G(t_i)/F^*(C_G(t_i))$ is an
elementary abelian group of order $\leqslant 4$.} \label{table}
\end{table}

\medskip

\item[(3)] Our next task is to make sure that $i$ is an involution of type $t_1$, $t_{n-1}$
 or $t_n$.
See Table~\ref{table} for the list of conjugacy classes of
involutions and their centralisers (the details are mostly taken
from \cite{lyons}).

For this purpose, construct  (almost) uniformly distributed random
elements of $C_G(i)$ as explained in \cite[Section~3]{old}. The
possibilities for $C_G(i)$ are given in Table~\ref{table}. Notice
that the groups $C_G(t_1)$, $C_G(t_{n-1})$ (in the group $G =
\Omega_{2n+1}(q)$) and $C_G(t_n)$ have much bigger orders of
semisimple elements than centralisers of involutions of other
types. Namely, in these three cases we can find in $C_G(i)$, with
probability $> O(1/n)$,  a `big' element, that is, an element of
order divisible by  a prime primitive divisor of $q^{n-1}\pm 1$.
However `big' elements are not present in the centralisers of
involutions of other types. As soon as a `big' element in
$C_G(i)$ is found,  we conclude that $i$ is of type $t_1$,
$t_{n-1}$  (in the group $G = \Omega_{2n+1}(q)$) or $t_n$ and move
to Step 4. If we cannot find a `big' element, we return to Step 1
and start our search for an involution of type $t_1$, $t_{n-1}$
or $t_n$ from the beginning.

\item[(4)] Now we want to make sure that we found an involution of
type $t_n$. For that we compute ${\rm pdrank}(ii^g)$ for $20n$
random elements $g$. We know from Lemmata~\ref{ortho},
\ref{lm:t1-in-PSp}  and \ref{lm:t1-tn-1} that if $i$ is not of
type $t_n$, or if $i$ is of type $t_n$ and $G = \Omega_{2n+1}(q)$
then ${\rm pdrank}(ii^g) \leqslant 6$. Hence if we found $g\in G$
such that ${\rm pdrank}(ii^g) > 6$ then we have a definite
witness against $G$ being the orthogonal group. In this case we
return ` $G$ is ${\rm Sp}_{2n}(q)$'.

If we have not found an element $g$ with ${\rm pdrank}(ii^g) >
6$, we can return to the search from Step 1.

\item[(4bis)] {\sc Shortcuts.} However, there are  some shortcuts which give a
quicker, although not one-sided solution.

\begin{itemize}
\item[(4i)] It follows from Lemma \ref{lm:t1-tn-1} that if ${\rm pdrank}(ii^g)
\leqslant 6$ for $O(n)$ random elements $g\in G$ and, for some
$g$, ${\rm pdrank}(ii^g)
> 4$ then $i$ is likely to be an involution of type $t_{n-1}$ in
the orthogonal group $\Omega_{2n+1}(q)$.

\item[(4ii)] It follows from Lemmata \ref{lm:t1-in-PSp} and \ref{lm:t1-tn-1}
 that if ${\rm pdrank}(ii^g)
\leqslant 4$ for  $O(n)$ random elements $g \in G$ and, for some
$g$, ${\rm pdrank}(ii^g) > 2$ then $i$ is likely to be an
involution of type $t_1$ in the orthogonal group
$\Omega_{2n+1}(q)$.

\item[(4iii)] It follows from follows from Lemmata \ref{lm:t1-in-PSp} and \ref{lm:t1-tn-1}
 that if\linebreak ${\rm pdrank}(ii^g)
\leqslant 2$ for  $O(n)$ random elements $g \in G$, then $i$ is
likely to be an involution of type $t_1$ in ${\rm Sp}_{2n}(q)$ or
of type $t_n$ in $\Omega_{2n+1}(q)$. Since the centralisers of
these involutions have very different structure (see
Table~\ref{table}), a further distinction between the two cases
can be made by running a black box for $C=C_G(i)$. The following
steps allow to detect a normal subgroup ${\rm SL}_2(q)$ in
$C=C_{{\rm Sp}_{2n}(q)}(t_1)$.

\begin{itemize}

\item    Choose $\varepsilon = \pm
1$ so that  $q^{n-1}-\varepsilon \equiv 2 \bmod 4$. Randomly
search in $C$ for an element $x$ of order divisible by $4$ and by
a primitive prime divisor of $q^{n-1}-\varepsilon$. If such an
element cannot be found after $O(n)$ attempts, $i$ is likely to
be an involution of type $t_n$ in $\Omega_{2n+1}(q)$.

\item Once the element $x$ is found,  compute an element $y$ of order
$4$ in $\langle x \rangle$.

\item Compute $O(n)$ elements $yy^c$ for random $c \in C$. If, for all of them,
 $o(yy^c) \mid q(q\pm 1)$ then $i$ is likely to
be an involution of type $t_1$ in ${\rm PSp}_{2n}(q)$ (and $y$
belongs to the ${\rm SL}_2(q)$-component of $C_G(i)$).
\item If, however, $o(yy^c) \nmid q(q\pm 1)$  for some $c \in C$
then $i$ is likely to be an involution of type $t_n$ in
$\Omega_{2n+1}(q)$.

\end{itemize}

\end{itemize}

\item[(5)] If we returned to Step 1 more than $20n$ times and have
not found a witness against the isomorphism $G \simeq
\Omega_{2n+1}(q)$, we conclude that $G \simeq \Omega_{2n+1}(q)$.

\end{itemize}

\section{Tori and semisimple elements in groups $B_n$ and $C_n$}

Definitions and notation used in the paper are standard and
borrowed from Carter \cite{car2,Carter-bk} and Gorenstein, Lyons
and Solomon \cite{lyons}.  We retain the notation and use the
results of \cite{old}. We use the description of maximal tori in
groups $B_n$ and $C_n$ as given in \cite[Section~2.1]{old}.

\subsection{Tori of order $\frac{1}{2}(q^n \pm 1)$.}

The crucial steps  of our  algorithm are aimed at finding
involutions belonging to
tori $\frac{1}{2}T_{q^n \pm 1}$ twisted by an element of the
cycle type $\pm n$. We have to distinguish the following cases:
\begin{itemize}
\item[$(i)$] $q^n \equiv -1 \bmod 4$, i.e. $n$ odd and $q \equiv -1  \bmod 4$.
\item[$(ii)$] $q^n \equiv 1  \bmod 4$, i.e. $n$ even or $q \equiv 1 \bmod 4$.
\end{itemize}
In the first case the maximal tori of order
$\frac{1}{2}(q^n  +1)$ contain involutions, in the second case the maximal tori
of order $\frac{1}{2}(q^n -1)$ contain involutions.

 We call tori $\frac{1}{2}T_{q^n\pm 1}$  {\it maximal twisted tori}.

We are looking for involutions belonging to maximal twisted tori.
For that purpose, we introduce the concept of {\it
$2$-height}\index{two@2-height}. Let $X$ be a group and $x \in X$
such that $o(x) = 2^km$ where $m$ is odd. Then $k$ is called the
{\it $2$-height of\/ $x$}. The $2$-height of $X$ is the maximum
of heights of elements in $X$.

\begin{lemma}
If $T$ is a torus of type $\frac{1}{2}T_{q^n\pm 1}$  then
 $$
 |N_G(T): T| = 2n.
$$
\label{lem:index}
\end{lemma}

\paragraph{{\it Proof.}} See \cite[Lemma~2.3]{old}.
\hfill $\square$

\subsection{Primitive prime divisors.}
\label{subsec:ppd}

Let $p$ be a prime number  $p$  and $m \geqslant 2$.  A prime
dividing $p^m-1$ but not $p^i-1$ for $1 \leqslant i < m$ is
called a {\it primitive prime divisor of $p^m-1$}. In our case
primitive prime divisors always exist by Zsigmondy \cite{zs}.
Furthermore an integer $j$ is called a $ppd(p;m)$ if $j\mid
q^m-1$ and $j$ is divisible by a primitive prime divisor of
$p^m-1$. Let $q = p^e$. A primitive prime divisor or $p^{em}-1$
will be called, for brevity, a $ppd(q,m)$. An element of $G$ of
order dividing a $ppd(q,n)$ is called a {\it primitive prime
divisor element}.

The role of primitive prime divisors becomes apparent from the
following simple lemma.

\begin{lemma}

\begin{itemize}
\item[{\rm (i)}] Let $r$ be a $ppd(q,n)$.
If  $r \mid q^m + 1$ for $m < n$, then $n$ is even and $m = n/2$.

\item[{\rm (ii)}]  Let $r$ be a $ppd(q,2n)$. Then $r \nmid q^m\pm 1$ for all
$m < n$ and $r \nmid q^n -1$.

\end{itemize}
\label{lm:ppd(q)}
\end{lemma}

\paragraph{{\it Proof.}} (i) Assume that $r \mid q^m + 1$ for some $m < n$. This
means that $r \mid q^{2n} -1$ and the Galois field
$\mathbb{F}_{q^{2m}}$ has a multiplicative subgroup $R$ of order
$R$. But the Galois field $\mathbb{F}_{q^n}$ also contains a
multiplicative subgroup of order $r$. Since the algebraic closure
$\overline{\mathbb{F}}_q$ contains a unique field of every give
order, the intersection $\mathbb{F}_{q^{2m}} \cap
\mathbb{F}_{q^n}$ is well defined. The inclusion
$\mathbb{F}_{q^n} < \mathbb{F}_{q^{2m}}$ because this would mean $
n \mid 2m$, contradicting  our assumption $m < n$. The equality
$\mathbb{F}_{q^{2m}} = \mathbb{F}_{q^n}$ leads to the desired
conclusion $ m = n/2$. Finally, we are left with the case when
$\mathbb{F}_{q^k} = \mathbb{F}_{q^{2m}} \cap \mathbb{F}_{q^n}$ is
a proper subfield of $\mathbb{F}_{q^n}$. But then the
multiplicative group of  $\mathbb{F}_{q^k}$ contains $R$, an
unique subgroup of order $r$ in the both $\mathbb{F}_{q^{2m}}$ and
$\mathbb{F}_{q^n}$, and therefore $r \mid q^k -1$, a
contradiction.

(ii) follows from (i); we need only notice that $r$ is odd, and,
since $q^n-1$ and $q^n+1$ have no odd  divisors in common, $r
\nmid q^n-1$ \hfill $\square$

\subsection{Good elements.}
\label{subsec:good}

Let $P$ be a maximal twisted torus of even order and let
$|P|=2^km$ where $k$ and $m$ are natural numbers
 and $m$ is odd.  We say that an element $g \in P$ is of {\it
  maximal 2-height}, if the 2-height of $g$ is $k$.  Since $P$ is
cyclic, elements of maximal 2-height exists. Furthermore half of
the elements of $P$ have maximal $2$-height as the cyclic group
generated by a 2-element $x$ of maximal $2$-height contains
$2^{k-1}$ elements of maximal $2$-height, namely $x^s$ for $s$ an
odd integer.

In the following we are looking for elements of maximal
$2$-height in a maximal twisted torus.
 We call $g$ {\it good} if it satisfies the
following conditions:
\begin{itemize}
\item ${\rm pdrank}(g) = 2n$
when $q^n \equiv -1 \bmod 4$ and ${\rm pdrank}(g) = n$ when $q^n
\equiv 1 \bmod 4$.

\item $g$ has maximal 2-height.

\end{itemize}

We shall call good elements belonging to maximal twisted tori in
$G$ {\em good twisted elements}.

\begin{lemma}
\label{shares} The number of good twisted elements in $G$ is at
least\/ $\frac{1}{5n}|G|$.
\end{lemma}

\paragraph{{\it Proof.}} Let $T$ be a maximal twisted torus in
$G$.
By Lemma~\ref{lem:index}
 $|N_G(T)/T| = 2n$.  Therefore there are $\frac{|G|}{2n|T|}$ tori conjugate to $T$.
Let $r$ be a $ppd(q,n)$ or $ppd(q,2n)$, depending on whether $T$
is of type $\frac{1}{2}T_{q^n- 1}$ or $\frac{1}{2}T_{q^n+ 1}$,
and denote by $T^{2r}$ the unique subgroup of index $2r$ in $T$.
Then any element whose image in $T/T^{2r}$ generate $T/T^{2r}$ is
a good element. It is easy to see that the number of such
elements in $T$ is $\frac{r-1}{2r}|T|$. Also, it follows from
Lemma~\ref{lm:ppd(q)} that a good element does not belong to a
proper subtorus of $T$ and thus belongs to a unique twisted
torus. Now we can estimate    the total number of good elements in
all twisted tori as at least
$$
\frac{|G|}{2n|T|} \cdot \frac{r-1}{2r}|T| = \frac{1}{4n}\cdot
\frac{r-1}{r}|G| \geqslant \frac{1}{5n} |G|,
$$
since, obviously, $r \geqslant 5$. \hfill $\square$

\subsection{Regular twisted elements.}
\label{subsec:regtwist}

A regular element in a torus of type $\frac{1}{2}T_{q^n\pm 1}$ or
will be called a {\em regular twisted\/} element. Notice that
good elements in these tori are regular twisted elements.

\begin{lemma}
\label{lem:share-reg-twist} The number of regular twisted
elements in $G$ conjugate to an element in a torus $T$ of the
given type $\frac{1}{2}T_{q^n-1}$ or  $\frac{1}{2}T_{q^n+1}$ is
at least $\frac{1}{4n}|G|$.
\end{lemma}

\paragraph{{\it Proof:}} See \cite[Lemma~2.5]{old}. \hfill
$\square$

\subsection{Involutions of type $t_n$ in ${\rm PSp}_{2n}(q)$.}
\label{subsec:tn}

For the following compare the description of conjugacy classes of
involutions in ${\rm
PSp}_{2n}(k)$ and ${\rm SO}_{2n+1}(k)$ as given in
Theorems~11.52 and 11.53 of Taylor \cite{tay} and Section~4.5 of \cite{lyons}.

The group ${\rm PSp}_{2n}(q)$ has $[\frac{n}{2}]+1$
conjugacy classes of involutions.
Let $s \in {\rm Sp}_{2n}(q)$ represent an element of order
2 in ${\rm PSp}_{2n}(q)$.
Then either $s^2 = -I$ or $s^2 = I$.  All involutions which are represented
by an element $s$ such that $s^2 = -I$ form one conjugacy class of
involutions of type $t_n$.  Furthermore for every even number
$2 \leqslant m \leqslant n$ there is one conjugacy class of
involutions of type $t_{m/2}$. If
$s$ is an involution of type $t_{m/2}$, then the eigenvalues of $s$
on the underlying natural module are $-1$ with multiplicity $m$ and $1$
with
multiplicity $2n-m$.

\begin{lemma}
\label{sympto} Let\/ $G:={\rm PSp}_{2n}(q)$ and\/ $i \in G$ an
involution of type $t_n$. If\/ $H$ is a  torus of\/ $G$ then
there exists an involution $j \in i^G$ such that\/ $H$ is
inverted by $j$.
\end{lemma}

\paragraph{{\it Proof.}} See \cite[Lemma~2.6]{old}.
\hfill $\square$

\begin{lemma}
\label{sympl} Let\/ $G:={\rm PSp}_{2n}(q)$ and\/ $i \in G$ be
an involution of type $t_n$. Then $i^gi^h$ is a regular twisted
element of odd order with probability $p \geqslant \frac{1}{20
n}$ for random $g, h \in G$.
\end{lemma}

\paragraph{{\it Proof.}}   See \cite[Lemma~2.7]{old}.
\hfill $\square$

\medskip

As an immediate corollary we have the following lemma.

\begin{lemma}
If $i$ is an involution of type $t_n$ in $G = {\rm PSp}_{2n}(q)$
then the share of elements $g\in G$ such that the element $ii^g$
is of odd order is at least $\frac{1}{20n}$.
\label{lem:odd-in-Sp}
\end{lemma}

\paragraph{{\it Proof:}} See \cite[Lemma~2.8]{old}.
\hfill $\square$

\subsection{Involutions of type $t_n$ in $\Omega_{2n+1}(q)$.}
\label{subsec:tn-in-O}

The group $\Omega_{2n+1}(q)$ has $n$ conjugacy classes of involutions.
There is one conjugacy class of involutions of type $t_{m/2}$ for each
even number $2 \leqslant m \leqslant 2n$. An involution $s$ is of
type $t_m$ if
the eigenvalues of $s$ on the underlying natural module are $-1$ with
multiplicity $m$ and $1$ with multiplicity $2n+1-m$.

\begin{lemma}
\label{ortho} Let\/ $G:=\Omega_{2n+1}(q)$ and\/ $i \in G$ an
involution of type $t_n$. Then  the element $i^g\cdot i^h$ acts
trivially on a subspace of codimension $2$ in the natural module
for $G$. In particular,  $o(i^gi^h) \mid q(q+1)$ or\/
$o(i^gi^h)\mid q(q-1)$ for all\/ $g, h \in G$.
\end{lemma}

\paragraph{{\it Proof.}} Let $V$ be the natural module for $G$.
Since $i$ is of type $t_n$ the dimension of the eigenspace $V_-$
of $i$ for the eigenvalue $-1$ is $2n$, while the eigenspace
$V_+$ for the eigenvalue $+1$ is $1$-dimensional. Therefore
$i^gi^h$ belongs to $C= C_{{\rm SL}(V)}(V_-g \cap V_-h)$ and $W =
V_-g \cap V_-h$ has codimension $\leqslant 2$ in $V$. Notice that
$C/O_p(C)$ is a subgroup of ${\rm SL}_2(q)$ (here $p = {\rm
char}\, \ff_q$). Notice also that the element $i^gi^h$ leaves
invariant the subspace $ V_+g + V_+h$. The result now easily
follows. \hfill $\square$ \vspace{1ex}

\begin{lemma}
Let\/ $G:=\Omega_{2n+1}(q)$ and\/ $i \in G$ an involution of type
$t_n$. Then the share of elements $g \in G$ such that the element $ii^g$ is
of odd order is at least
$$
\frac{1}{4(1+ 2/q +o(1/q))}.
$$
\label{lem:odd-in-O}
\end{lemma}

\paragraph{{\it Proof:}} We use notation of the previous lemma with $W = V_- \cap V_-g$
and $U = V_+  + V_+g$. Let $q \equiv 1 \bmod 4$. Then a direct computation with the help of
information about finite orthogonal geometries contained in Section III.6 of  Artin \cite{artin}
shows that, with  probability
$$
\rho > \frac{1}{4(1+ 2/q +o(1/q))},
$$
the subspace $U$ is nondegenerate and has index $0$.
Hence the group $\langle i, i^g\rangle $
is isomorphically imbedded into $O(U) \simeq \zz_{q+1} \rtimes \zz_2$
and half of the elements $ii^g$ have odd order, yielding the desired estimate.
The case $q \equiv -1 \bmod 4$ can be treated
analogously. \hfill $\square$

\subsection{Involutions of type $t_1$ in ${\rm PSp}_{2n}(q)$}

\begin{lemma} If\/ $i$ is an involution of type $t_1$ in\/ ${\rm
PSp}_{2n}(q)$ then
$$o(i^gi^h) \mid q(q \pm 1).$$

\label{lm:t1-in-PSp}
\end{lemma}

\paragraph{{\it Proof:}}
Similarly to the previous section, a preimage in $G={\rm
Sp}_{2n}(q)$ of the involution $i$  (which we denote by the same
symbol $i$) is a `low dimensional' involution, that is, the
eigenspace $V_-$ has dimension $2$ or $n-2$. Assume that $\dim
V_- = 2$, the other case is similar. Then, for random $g,h \in G$,
the group $\langle i^g, i^h\rangle$ leaves the subspaces $V_-g +
V_-h$ and $V_+g \cap V_+h$ invariant and $V = (V_-g + V_-h)
\oplus (V_+g \cap V_+h)$. Denote $W = V_-g + V_-h$. Notice that
$\langle i^g, i^h\rangle < {\rm Sp}(W)$.

We have three possibilities: $\dim W = 2$, $3$ or $4$. If $\dim W
= 2$ the $i^g = i^h$. If the restriction of the symplectic form
to $W$ is degenerate (this automatically holds when $\dim W = 3$)
then ${\rm Sp}(W)/O_p({\rm Sp}(W))\simeq {\rm Sp}_2(q) = {\rm
SL}_2(q)$, and the lemma follows. If $\dim W = 4$ and the
restriction of the symplectic form to $W$ is non-degenerate then
${\rm Sp}(W)= {\rm Sp}_4(q)$ and $i^g$ and $i^h$ are involutions
of type $t_1$ in ${\rm Sp}(W)$. Now the result follows from an
easy computation in ${\rm Sp}_4(q)$. \hfill $\square$

Moreover, the same analysis yields the following result.

\begin{lemma}
If\/ $i$ is an involution of type $t_1$ in\/ ${\rm PSp}_{2n}(q)$
then the share of elements in $G$ such that\/ $ii^g$ is of odd
order is at least\/ $1/4$. \label{lm:t1-in-Psp-odd}
\end{lemma}

\subsection{Involutions of type $t_1$ and $t_{n-1}$ in
$\Omega_{2n+1}(q)$}

\begin{lemma}
Let\/ $i$ be an involution of type $t_1$ or\/ $t_{n-1}$ in
$G=\Omega_{2n+1}(q)$ and\/ $g$ is a random element in $G$.
\begin{itemize}
\item[{\rm (i)}] If\/ $i$ is of type $t_1$ then\/ ${\rm pdrank}(ii^g) \leqslant 4$.
Moreover, ${\rm pdrank}(ii^g) > 2 $ with probability bounded from
below by a constant.
\item[{\rm (ii)}] If\/ $i$ is of type $t_{n-1}$ then\/ ${\rm pdrank}(ii^g) \leqslant 6$.
Moreover, ${\rm pdrank}(ii^g) > 2$ with probability bounded from
below by a constant.

\item[{\rm (iii)}] If\/ $i$ is of type $t_1$ or $t_{n-1}$ then the element $ii^g$ has
odd order with probability bounded from below by a constant.
\end{itemize}
\label{lm:t1-tn-1}
\end{lemma}

The experimental evidence suggests that the constants in (i) and (ii) are $1/8$
for $t_1$ and $1/4$ for $t_{n-1}$.

\paragraph{{\it Proof.}} The proof is similar to the previous
result. \hfill $\square$

\subsection{Involutions produced by good elements.}

\begin{lemma}
\label{type} Let\/ $i$ be the involution in a maximal twisted
torus\/ $T < G$ . Then\/ $i$ is of type $t_n$.

\end{lemma}

\paragraph{{\it Proof.}} See \cite[Lemma~2.13(1)]{old}.   \hfill $\square$

We shall use the following corollary from this result and
Lemma~\ref{shares}.

\begin{lemma}
Let\/ $g$ be a random element in $G$. Then the involution\/ $i$
produced by $g$ has type $t_n$ with probability at\/ least\/
$1/5n$.
\end{lemma}

Finally,

\section{Detailed description of the algorithm}

\subsection{Step 1: search for a good element $x$.}
The algorithm spends most of its timing testing if elements are
good.  This is done exactly as described in
\cite[Section~4]{old}.

It easily follows from Lemma~\ref{shares} that a sample of $N$
random elements of $G$
 contains a good element with probability
 $$
 1-\left(1-\frac{1}{5n}\right)^{N}  > 1- e^{-\frac{N}{5n}}.
 $$
 Therefore if we are given fixed margin of error $\epsilon$,
$0 < \epsilon < 1$, and take a sample of size $O(n\log
(1/\epsilon))$, then we can be satisfied that it contains a good
element with probability $1-\epsilon$. As proven in \cite{old},
production of a random element and testing it for goodness require
polynomial time.

\subsection{Step 2: involution $i = {\rm i}(x)$.}
As soon as we have a good element $x$ and know that it is an
element of maximal $2$-height in a cyclic group of order
$\frac{1}{2}(q^n\pm 1)$, the involution $x={\rm i}(x)$ can be
obtained by raising $x$ to the power $\frac{1}{4}(q^n\pm1)$. As
discussed in \cite{old}, this requires  polynomial time.
Lemma~\ref{type} asserts that $i$ is  of  type $t_n$.

\subsection{Step 3: involution of type $t_n$.}

Notice that if $n$ is odd then the involution $i$ is of type $t_n$.

However, if $n$ is even then $i$ can happen to be of different
type.

We can check whether $i$ is of type $t_n$ or not by constructing
uniformly distributed random elements   $\zeta_1(x)$ of  $C_G(i)$
as described in \cite[Section~3]{old}. If $i$ is of type $t_1$,
$t_{n-1}$ (in $\Omega_{2n}^+(q)$) or $t_n$ then, in view of
Lemmata~\ref{lm:t1-in-Psp-odd}, \ref{lm:t1-tn-1} and
\ref{lem:odd-in-O}, this method produces (in polynomial of $n\log
q$ time) random and uniformly distributed elements of $C_G(i)$.
Let $L = F*(C_G(i))$. Depending on whether $i$ is one of the
$t_1$, $t_{n-1}$ (in $\Omega_{2n}^+(q)$) or $t_n$, $L$ contains
'big' elements, that is, elements of order divisible by a
primitive prime divisor of $q^{n-1}-1$, $q^{n-1}+1$, $q^{n-1}+1$
(compare Table~\ref{table} and Table~\ref{table:orders}). It can
be observed directly from the formulae for orders in
Table~\ref{table:orders} that centralisers of involutions of
other types do not contain `big' elements.

\begin{table}[h]
\begin{tabular}{|| c | c ||}
\hline & \\[1ex]
$|{\rm GL}_k(q)|$ & $q^{k(k-1)/2}(q-1)(q^2-1)\cdots (q^k-1)$ \\[2ex]
\hline & \\[1ex]
 $|{\rm PSp}_{2k}(q)| = |\Omega_{2k+1}(q)|$ & $\frac{1}{(2,q-1)}
q^{k^2}(q^2-1)(q^4-1)\cdots (q^{2k}-1)$ \\[2ex]
\hline & \\[1ex]
$|{\rm P}\Omega^+_{2k}(q)|$ & $\frac{1}{(4,q^k-1)}
q^{k(k-1)}(q^2-1)(q^4-1)\cdots (q^{2k-2}-1)(q^k-1)$ \\[2ex]
\hline & \\[1ex]
$|{\rm P}\Omega^-_{2k}(q)|$ & $\frac{1}{(4,q^k+1)}
q^{k(k-1)}(q^2-1)(q^4-1)\cdots (q^{2k-2}-1)(q^k+1)$ \\[2ex]
\hline & \\[1ex]
$|{\rm GU}_k(q^2)|$ &
$q^{k(k-1)/2}(q+1)(q^2-1)(q^3+1)(q^4-1)\cdots
(q^k-(-1)^k)$\\[2ex] \hline

\end{tabular}

 \caption{Orders of some classical groups, \cite[pp.~2--8]{car2}.} \label{table:orders}
\end{table}

It is easy to see that `big'
elements are regular elements in maximal tori of orders
$\frac{1}{2}(q \pm 1)(q^{n-1} \pm 1)$ or $\frac{1}{2}(q^n +1)$.
Using the same technique as in Sections~2.1--2.3, we can see that
the probability of a random element to be `big'  is $O(1/n)$ and
hence one of them has to pop up from the black box for $C_G(i)$
after $O(n)$ iterations. As soon as we have found a `big' element,
 we know that $i$ is of type  $t_1$, $t_{n-1}$ (in $\Omega_{2n+1}(q)$)
  or $t_n$ and can go to Step 4.

If a 'big' element cannot be found after $O(n)$ tests of elements
from $C_G(i)$, then to Step 1 and start the search for an
involution of type $t_1$ or $t_n$ from the beginning.

\subsection{Step 4.}

Now we want to make sure that we found an involution of type
$t_n$ in ${\rm PSp}_{2n}(q)$. For that we compute ${\rm
pdrank}(ii^g)$ for $20n$ random elements $g$. We know from
Lemma~\ref{sympl} that if $i$ is of type $t_n$ in ${\rm
PSp}_{2n}(q)$, then, with probability $> 1-e^{-1}$, the primitive
divisor rank of one of the elements $ii^g$ is at least $2n$. We
know from Lemmata~\ref{ortho}, \ref{lm:t1-in-PSp} and
\ref{lm:t1-tn-1} that if $i$ is not of type $t_n$, or if $i$ is
of type $t_n$ and $G = \Omega_{2n+1}(q)$ then ${\rm pdrank}(ii^g)
\leqslant 6$. Hence if we found $g\in G$ such that ${\rm
pdrank}(ii^g) > 6$ then we have a definite witness against $G$
being the orthogonal group. In this case we return ` $G$ is ${\rm
Sp}_{2n}(q)$'.

If we have not found an element $g$ with ${\rm pdrank}(ii^g) >
6$, we can return to the search from Step 1.

The estimate of the number of elements $ii^g$ needed in this step
is, of course, very crude; an experimental evidence suggests
that, instead of $20n$ elements $ii^g$, $8$ elements suffice.

{\sc Shortcuts} in Step 4 are mostly self-explanatory.

\subsection{Step 5.}

The involution $i$ of type $t_n$ behaves very differently in ${\rm PSp}_{2n}(q)$
and $\Omega_{2n+1}(q)$.
Namely, in ${\rm PSp}_{2n}(q)$ an element $i^gi^h$ for random $g,h \in G$
 is a regular twisted element  with probability
$\geqslant \frac{1}{20n}$ (Lemma~\ref{sympl}), while in $\Omega_{2n+1}(q)$
we always have $o(i^gi^h) \mid q(q \pm 1)$ (Lemma \ref{ortho}).
Obviously, a regular twisted element cannot have order dividing $q(q\pm 1)$.
As soon as we have that
$o(i^gi^h)$ does not divide $q(q \pm 1)$,
the pair $(i^g, i^h)$ is a witness against the isomorphism
$G \simeq \Omega_{2n+1}(q)$ and in that case $G \simeq {\rm PSp}_{2n}(q)$.

On the other hand, if in  $20n$ tests we had the same result
$o(i^gi^h) \mid q(q \pm 1)$, then $G \simeq \Omega_{2n+1}(q)$
with probability $\geqslant 1 - e^{-1}$. Timing estimate here is
very crude and can be considerably improved. In practice,  if
$G\simeq {\rm PSp}_{2n}(q)$ then $o(i^gi^h) \nmid q(q \pm 1)$ at
the very first attempt.

\section{The GAP code.}

The algorithm was tested in GAP\footnote{The software package GAP is available,
for variety of platforms, from http://aldebaran.math.rwth-aachen.de/~GAP/}.
The Appendix contains a (very
primitive) GAP source code {\tt Bn-v-Cn.g}, which generates (at
random) a group $G={\rm Sp}_{2n}(q)$ or ${\rm SO}_{2n+1}(q)$ and
initiates the functions {\tt Test(G)} and  {\tt PlainTest(G)}
which determine the type of the group $G$ with the help of the
algorithm described in this paper, with ({\tt Test(G)}) or
without ({\tt PlainTest(G)}) shortcuts in Step 4.

The functions {\tt PlainTest(G)} and {\tt Test(G)} return the list of $4$ integers.
Of these,
\begin{itemize}
\item The first position is 0 or 1, according to the answer
    {\em symplectic} or *{\em orthogonal} given by the algorithm.
\item The second position is the rule number used in this decision
(used for the debugging of the algorithm).
\item The third position is the dimension of the eigenspace for the eigenvalue $-1$
of the involution used in the identification of the group
this information is used only for the debugging of the algorithm).
\item The fourth number  is the numbers of re-runs of the test
 if the first run was inconclusive.
 \end{itemize}

I wish to emphasise that this GAP code was written only for
checking the algorithm and is not designed for the practical use. In
particular, it works with the groups ${\rm Sp}_{2n}(q)$ and ${\rm
SO}_{2n+1}(q)$ instead of ${\rm PSp}_{2n}(q)$ and
$\Omega_{2n+1}(q)$, although all computations are done in the way
which makes sense in the factor group modulo the scalars and in
the commutator group; in particular, it works, instead of
involutions, with {\em pseudoinvolutions}, that is, elements of order
$2$ or $4$ whose images modulo the scalars are involutions.

To see how the code works, enter at the GAP command prompt:

\medskip
\noindent
 {\tt
gap> Read("Bn-v-Cn.g"); test:=Test(G); Interpret(test); }

\medskip
\noindent or

\medskip
\noindent
 {\tt
gap> Read("Bn-v-Cn.g"); test:=PlainTest(G); Interpret(test); }

\medskip

To check the answer, enter

\medskip
\noindent
 {\tt
gap>  G;  }

\medskip

I have to warn, however, that {\tt PlainTest(G)}
takes considerable time to produce the answer if $G$ is
the orthogonal group.

So far more than $2000$ experiments with the function {\tt
Test(G)} have not yield a single wrong answer,
while {\tt PlainTest(G)} erred once in $600$ experiments.

This was my first ever attempt to write in GAP, and I bravely
ignored the good programming practice; in particular, the global
variables q, n, type (0 or 1 for symplectic or orthogonal) for
generation of a black box group can be changed only in the file
{\tt Bn-v-Cn.g}, which after that has to be read again to generate
all the required functions.

Also, I am using the function {\tt Order(element)}, to make the code simpler;
this is not necessary, of course, and it can be replaced by
raising the element in the appropriate powers as discussed in \cite{old}.
This would make it possible to work in "large" groups.

Structurally, I was not concerned with the optimal branching in
the program. I intentionally wanted all branches being regularly
visited, although, in many cases, this could be avoided.

\vfill

\noindent
Department of Mathematics, UMIST,
PO Box 88, Manchester M60 1QD, United Kingdom\hfil\break
\textsc{E-mail:} \texttt{Alexandre.Borovik@umist.ac.uk}\hfil\break
\textsc{Web:} \texttt{http://www.ma.umist.ac.uk/avb/}
\hfil\break

\newpage
\section*{Appendix: GAP code}

\small
\begin{verbatim}

####################################################################
###
### Bn-v-Cn.g
###
### A. V. Borovik, April 2001
###
### The algorithm tests involutions in orthogonal and symplectic groups
### thus distinguishing between the two.
### Given a linear group G = PSp(2n,q) or Omega(2n+1, q),
### with known q and n,
### the algorithm determines whether G is sympelctic or
### orthogonal.
### The algorithm implemented here contains a fix for a bug
### in the original theoretic version,
### C. Altseimer and A. Borovik, "Probabilistic recognition
### of orthogonal and symplectic groups",
### Groups and Computations III, (eds W.~M.Kantor and A.
### Seress), de Gruyter, Berlin and New York, 2001, pp. 1-20.
###
### This code was not designed for the practical use,
### but just for checking the theoretical algorithm.
### In particular, since I had no means for generation of
### black boxes for PSp(2n,q) for large n and q,
### the algorithm actually works with Sp(2n,q) instead of
### PSp(2n,q), but make only those group calculations
### which make sense in its factor group PSp(2n,q).
### In particular, it works with "pseudoinvolutions",
### elements of order 2 or 4 which have order 2 in G/Z(G).
###
### Recommended use:
### gap> Read("Bn-v-Cn.g"); test:=Test(G); Interpret(test);
### or
### gap> Read("Bn-v-Cn.g"); test:=PlainTest(G);  Interpret(test);
### To check the answer, enter
### gap> G;
####################################################################


### Global constants: q, n, type, G, effort

q:= Random([5,7,9,11,13,17,19,23,25,27,31]);


if q >19 then               ### selection of the dimension of
    n:= Random([3..6]);     ### of the group with the aim to
fi;                         ### keep the exponent at
if q > 10 and q < 20 then   ### computationally feasible size
    n:= Random([3..12]);
fi;
if q < 10 then
    n:= Random([3..20]);
fi;



type:= Random([0,1]);   ### type = 0 means symplectic
                        ### type = 1 means orthogonal
if type = 1 then
    G:=SO(2*n+1,q);     ### Generation of the group G, type =1
else
    G:=Sp(2*n,q);       ### Generation of the group G, type = 0
fi;

effort := 2;            ### This parameter determines how hard
                        ### we are looking for the witness
                        ### elements. Most searches make
                        ### effort * n attempts to find a
                        ### "big" element
                        ### Effort =4 corresponds to the theoretical threshold:
                        ### in the group Sp(2n,q) or Omega(2n+1,q),
                        ### a random element is a regular element
                        ### in a maximal twisted torus of order
                        ### q^n -1 or q^n+1 with probability > 1/4n.

###
### Functions
###


OddPart := function(nn)
### Odd part of the integer nn
    local m;
    m :=nn;
    while m mod 2 = 0 do
        m := m/2;
    od;
    return m;
end;


EvenPart := function(nn);
### Even part of the integer nn

    return nn/OddPart(nn);
end;


PrimitiveDivisor := function(qq,nn)
### Finds the largest divisor of the number qq^nn -1
### which has no non-trivial divisors in common with
### any of the numbers qq^l - 1 for l < nn.

    local pd,ll,gcd,ii;
    pd := qq^nn-1;
    for ii in [1..nn-1] do
        ll := qq^ii -1;
        gcd := Gcd(pd,ll);
        while gcd > 1 do
            pd := pd/gcd;
            gcd:= Gcd(pd,gcd);
        od;
    od;
    return pd;
end;



PseudoInvolution := function(g)
### Finds an element in the cyclic group <g>
### which is an involution in G/Z(G)
### The function uses the standard GAP function Order(g)
### which uses prime factorisation and does not work in big
### groups. However, it is easy to eliminate calls to Order
### and replace them by exponentiation.

    local order, i, rank;
    order := Order(g);
    if order mod 2 = 1 then
        return g^0;
    else
        i := g^(order/2);
    fi;

    rank := Rank(i+i^2);
    if (rank = 0 and type = 0) then
        if order mod 4 <> 0 then
            return g^0;
        else i := g^(order/4);
        fi;
    fi;
    return i;
end;


InCentraliser := function(j)
### Returns an element from the centraliser C_G(j)
### Most of the elements produced by this function are
### involutions; their distribution in C_G(j) is
### invariant under action of C_G(j) by conjugation
### Some (relatively small) portion of elements is
### produced by a different algorithm $\zeta_1$ in
### notation of our paper; these elements are distributed
### uniformly in the group C_G(j).

    local gg,y,o,z,c;
    gg := PseudoRandom(G);
    y := j * j^gg;
    o := Order(y);
    if j^4 <> j^2 then
        if (o = 1 or o = 2 or o = 4) then
            z := y;
        fi;

        if o mod 4 = 0 then
            if y^(o/2) = j^2  then
                z := y^(o/4);
            else
                z:= y^(o/2);
            fi;
        fi;

        if o mod 4 = 2 then
            if y^(o/2) = j^2  then
                z := gg*y^(((o/2)-1)/2);
            else
                z:=j;
            fi;
        fi;

        if o mod 2 = 1 then
            z := gg*y^((o-1)/2);
        fi;

    else
        if o mod 2 = 0 then
            z := y^(o/2);
        else
            z := gg*y^((o-1)/2);
        fi;
    fi;
    return z;
end;


OccurrencesOfPrimitiveDivisors := function(m)
### Finds all integers k such that m is divisible by
### a primitive prime divisor of q^k-1.

    local tuple, k;
    tuple := 0*[1..2*n];
    for k in [1..2*n] do
        if Gcd(PrimitiveDivisor(q,k),m) > 1 then
            tuple[k] :=k;
        fi;
    od;
    return tuple;
end;

PrimitiveDivisorRankOfElement := function(g)
### Finds the maximal integer k =< 2*n such that
### the order of the element g is divisible
### by a primitive prime divisor of q^n-1
### (remember, q and n are global constants).
### Again, calls to Order(g) can be eliminated

    return Maximum(OccurrencesOfPrimitiveDivisors(Order(g)));
end;


IsGood := function(g)
### Tests whether the element g is good or not.

    if (q^n -1) mod 4 = 0 then
        return (PrimitiveDivisorRankOfElement(g) = n and EvenPart(Order(g)) =

(EvenPart(q^n-1)/(2^type)));
    fi;
    if (q^n + 1) mod 4 = 0 then
        return (PrimitiveDivisorRankOfElement(g) = 2*n and EvenPart(Order(g))=

(EvenPart(q^n+1)/(2^type)));
    fi;
end;

GoodElement := function(G)
### Finds a good element in the group G,
### of course, under the assumption that
### G = Sp(2*n,q) or S0(0,2*n+1,q);
### in the latter case, the good element
### is actually found in the commutator group
### Omega(2*n+1,q) of S0(0,2*n+1,q).

    local g1, g2, g;
    if type = 1 then
        g1:=PseudoRandom(G);
        g2:=PseudoRandom(G);
        g :=g1^-1*g2^-1*g1*g2;
        while IsGood(g) = false do
            g1:=PseudoRandom(G);
            g2:=PseudoRandom(G);
            g :=g1^-1*g2^-1*g1*g2;
        od;
    else
        g:=PseudoRandom(G);
        while IsGood(g) = false do
            g:=PseudoRandom(G);
        od;
    fi;
    return g;
end;


TypeOfCentraliser := function(j)
### Returns one of the lists [], [1], [2], and [1,2],
### depending on whether an element of order dividing
### a primitive prime divisor of q^(n-1)-1 or a
### primitive prime divisor of q^(n-1)+1 is found in C_G(j)
### (this adds 1 or 2 to the list, correspondingly).
###
### The procedure is pretty slow; it can be made faster
### by a better way of generating random elements in C_G(j)
### (for example, the use of the Product Replacement Algorithm
### should considerably speed it up).

    local x, maxppd, typeOfC, counter;
    x := InCentraliser(j);
    maxppd :=0;
    typeOfC := Set([]);
    counter := 0;
    while typeOfC <> Set([1,2]) and counter < effort*n do
        x := x * InCentraliser(j);
        maxppd := PrimitiveDivisorRankOfElement(x);
        if maxppd = n-1  then
            Add(typeOfC, 1);
        fi;
        if maxppd = 2*(n-1) then
            Add(typeOfC, 2);
        fi;
        counter := counter +1;
    od;
    return Set(typeOfC);
end;

ProductsOfConjugates := function(j)
### Determines how big is the primitive divisor rank
### of the product j * j^x of the element j and
### its random conjugate  j^x.

    local x, x2, y, maxppd, counter;
    maxppd := 0;
    counter := 0;
    if type = 0 then
        while (maxppd < 8 and counter < effort*n) do
            x := PseudoRandom(G);
            y := j * j^x;
            maxppd := Maximum([maxppd, PrimitiveDivisorRankOfElement(y)]);
            counter := counter +1;
        od;
    else
        while (maxppd < 8 and counter < effort*n) do
            x := PseudoRandom(G);
            x2 := PseudoRandom(G);
            x := x^-1 * x2^-1 * x * x2;
            y := j * j^x;
            maxppd := Maximum([maxppd, PrimitiveDivisorRankOfElement(y)]);
            counter := counter +1;
        od;
    fi;
    return maxppd;
end;


CriticalInvolution := function(G)
### Finds an involution j in the cyclic group generated by a good element,
### with the extra property that the centraliser C_G(j) has an element
### whose order is divisible by a primitive prime divisor of q^(n-1)-1
### or q^(n-1)+1.

    local x, j, typeOfC;
    typeOfC := [];
    x := GoodElement(G);
    j := PseudoInvolution(x);
    while (typeOfC = Set([]) or j=j^2)  do
        x := GoodElement(G);
        j := PseudoInvolution(x);
        typeOfC := TypeOfCentraliser(j);
    od;
    return j;
end;



PreTest := function(G)
### Tests the group G for being symplectic or orthogonal
### The result may be inconclusive, in which case
### the test has to be repeated.

local g, j, typeOfC, conj, conj2, x, y, z, minus, plus, maxppd,
evenpart, counter, maxppd2,counter2;

    ### The case when n =3

    if n  = 3 then
        g := GoodElement(G);
        j := PseudoInvolution(g);
        conj:= ProductsOfConjugates(j);
        counter :=0;
        while conj < 3 and counter < effort * n do
            conj2 := ProductsOfConjugates(j);
            conj := Maximum([conj, conj2]);
            counter := counter + 1;
        od;
        if conj < 3 then
            return [1, 13, Rank(j-j^2)/2];
        else
            return [0, 13, Rank(j-j^2)/2];
        fi;
    fi;

    ###     The general case
    j := CriticalInvolution(G);
    conj := ProductsOfConjugates(j);
    if  conj > 7 then
        return [0, 1, Rank(j-j^2)/2];
    fi;

    typeOfC := TypeOfCentraliser(j);
    if typeOfC = [1,2] and conj = 4 then
        return [1,2, Rank(j-j^2)/2];
    fi;
    if n mod 4 = 2 and q mod 4 =1 and typeOfC = [1] and conj = 6 then
        return [1,3, Rank(j-j^2)/2];
    fi;
    if n mod 4 = 2 and q mod 4 =3 and typeOfC = [2] and conj = 6 then
        return [1,4, Rank(j-j^2)/2];
    fi;
    if n mod 4 = 0 and typeOfC = [1,2] and conj < 3 then
        return [1,5, Rank(j-j^2)/2];
    fi;
    if n mod 2 = 1 and typeOfC = [1,2] and conj < 3 then
           return [1,6, Rank(j-j^2)/2];
    fi;

    if  n mod 2 = 0 and typeOfC = [1,2] and conj < 3 then
        ####
        minus := (EvenPart(q^(n-1)-1) mod 4 <> 0);
        plus := (EvenPart(q^(n-1)+1) mod 4 <> 0);
        x := InCentraliser(j);

        if minus then
            maxppd :=0;
            evenpart :=1;
            counter:= 0;
            while (maxppd <> n-1 or evenpart < 4) and counter < 3*effort*n do
                x := x * InCentraliser(j);
                maxppd := PrimitiveDivisorRankOfElement(x);
                evenpart := EvenPart(Order(x));
                counter := counter +1;
            od;
            if (maxppd <> n-1 or evenpart < 4) then
                return [1,7, Rank(j-j^2)/2];
            fi;

        fi;

        if plus then
            maxppd :=0;
            evenpart :=1;
            counter:= 0;
            while (maxppd <> 2*(n-1) or evenpart < 4) and counter < 3*effort*n  do
                x := x * InCentraliser(j);
                maxppd := PrimitiveDivisorRankOfElement(x);
                evenpart := EvenPart(Order(x));
                counter := counter +1;
            od;
            if (maxppd <> 2*(n-1) or evenpart < 4) then
                return [1,8, Rank(j-j^2)/2];
            fi;
        fi;

        y := x^(Order(x)/evenpart);
        counter2:=0;
        maxppd2:=0;
        while (maxppd2 < 3 and counter2 < effort*n) do
            x := x * InCentraliser(j);
            z := y * y^x;
            maxppd2 := Maximum([maxppd2, PrimitiveDivisorRankOfElement(z)]);
            counter2 := counter2+1;
        od;

        if maxppd2 < 3 then
            return [0,9, Rank(j-j^2)/2];
        fi;
        if maxppd2 > 2 then
            return [1,10, Rank(j-j^2)/2];
        fi;

        ####
    fi;

    if n > 3 then
        return [-1,11, Rank(j-j^2)/2];
    fi;

end;

Test := function(G)
### PreTest(G) is repeated  until we get a conclusive answer

    local pretest, counter;
    pretest := PreTest(G);

    if pretest[1] = 0 or  pretest[1] = 1 then
        Add(pretest, 0);
        return pretest;
    fi;

    counter :=0;
    while pretest[1] = -1 do
        pretest := PreTest(G);
        counter := counter + 1;
    od;
    #Add(pretest, 0); Add(pretest, 0); Add(pretest, 0);
    Add(pretest, counter);
    return pretest;
end;



PlainPreTest := function(G)

### Tests the group G for being
### symplectic or orthogonal.
### No shortcuts in the algorithm.
### The result may be inconclusive,
### in which case the test has to be repeated.

    local g, j, typeOfC, conj;
    g := GoodElement(G);
    #g := PseudoRandom(G);
    j := PseudoInvolution(g);
    while j = j^0 do
        g := GoodElement(G);
        #g := PseudoRandom(G);
        j := PseudoInvolution(g);
    od;
    conj := ProductsOfConjugates(j);
    typeOfC := TypeOfCentraliser(j);
    if  (conj > 7 and typeOfC <> []) or (conj > 2 and n = 3) then
            return [0, 14, Rank(j-j^2)/2];
    else
            return [-1, 14, Rank(j-j^2)/2];
    fi;
end;

PlainTest := function(G)

### PlainPreTest(G) is repeated until we get a conclusive answer.

    local pretest, counter;
    pretest := PlainPreTest(G);

    if pretest[1] = 0 or  pretest[1] = 1 then
        Add(pretest,0);
        return pretest;
    fi;

    counter :=0;
    while pretest[1] = -1 and counter < effort*n do
        pretest := PlainPreTest(G);
        counter := counter + 1;
    od;
    #Add(pretest, 0); Add(pretest, 0); Add(pretest, 0);
    Add(pretest, counter);
    if counter > effort*n-1 then
        pretest[1] := 1;
    fi;
    return pretest;
end;

Interpret := function(test)
### Interpretation of the test result

    local l;
    l := Length(test);
    if test[1] = 0 then
        Print(" The group is symplectic.\n");
    fi;

    if test[1] =1 then
        Print(" The group is orthogonal.\n");
    fi;

    if test[2] = 1 then
        Print(" This is the DEFINITE answer because the test has found an involution j \n");
        Print(" whose centraliser  C_G(j)  contains an element of\n");
        Print(" order divisible by a primitive prime divisor of q^(n-1)-1  \n");
        Print(" or  q^(n-1)+1, and such that, in addition, there is a product\n");
        Print(" j*j^x  of two involutions conjugate to  j  whose order is divisible\n");
        Print(" by a primitive prime divisor of  q^k-1 for some  k>6;\n");
        Print(" this is possible only for an involution of type t_n in the symplectic group.\n");
    fi;

    if test[2] = 2 then
        Print(" An involution  j was found whose centraliser  C_G(j) contains\n");
        Print(" elements of order divisible by primitive prime divisors of the both\n");
        Print(" q^(n-1)-1  or  q^(n-1)+1, and with the maximum of the primitive prime\n");
        Print(" divisor ranks of random products of conjugates j^x * j^y being \n");
        Print(" equal 4. This almost certainly identifies this involution j \n");
        Print(" as involution of type t_1 in SO(2n+1,q).\n");
        Print(" To be on the safe side, repeat the test.\n");
    fi;

    if test[2] = 3 or test[2] = 4 then
        Print(" In short, the test found an involution which behaves as an\n");
        Print(" involution of type t_{n-1} in SO(2n+1,q).\n");
        Print(" this is not a very conclusive answer; to be on the safe side,\n");
        Print(" repeat the test.\n");
    fi;

    if test[2] = 5 or test[2] = 6 then
        Print(" In short, the test found an involution which behaves as an\n");
        Print(" involution of type t_n in SO(2n+1,q).\n");
        Print(" This is not a very conclusive answer; to be on the safe side,\n");
        Print(" repeat the test.\n");
    fi;

    if test[2] = 7 or test[2] = 8 or test[2] = 10 then
        Print(" The test found an involution j which behaves either as an involution\n");
        Print(" of type t_1 in PSp(2n,q), or as an involution of type t_n in\n");
        Print(" Omega(2n+1,q); further distinction was made by not being able to find\n");
        Print(" a 2-element from the component of type SL(2,q) which should be present\n");
        Print(" in the centraliser C_G(j) in the sympelctic case.\n");
        Print(" To be on the safe side, repeat the test.\n");
    fi;

    if test[2] = 9 then
        Print(" The test found an involution j which behaves either as an involution\n");
        Print(" of type t_1 = t_{n-1} in PSp(2n,q), or as an involution of type t_n in\n");
        Print(" Omega(2n+1,q); further distinction was made by being able to find\n");
        Print(" a 2-element in C_G(j) which behaves as if it belongs to\n");
        Print(" the component of type SL(2,q) which should be present\n");
        Print(" in the centraliser C_G(j) in the symplectic case.\n");
        Print(" To be on the safe side, repeat the test.\n");
    fi;

    if test[1] = 0 and test[2] = 13 then
        Print(" Since n = 3, the test is much simpler: we found an involution j\n");
        Print(" generated by a good element and such that the primitive prime divisor rank\n");
        Print(" of some random product  j^x * j^y of conjugates of j is at least 3. \n");
        Print(" This DEFINITELY identifies the group as symplectic.\n");
    fi;

    if test[1] = 1 and test[2] = 13 then
        Print(" Since n = 3, the test is much simpler: we found an involution j\n");
        Print(" generated by a good element with the  primitive prime divisor ranks\n");
        Print(" of a random products of conjugates j^x * j^y at most 2. \n");
        Print(" This almost definitely identifies the group as orthogonal.\n");
        Print(" To be on the safe side, repeat the test.\n");
    fi;

    if test[1] = 0 and test[2] = 14 then
        Print(" The programme found an involution which behaves as \n");
        Print(" involution  of type t_n in PSp(2n,q). \n");
        Print(" This is a DEFINITE answer given by  \n");
        Print(" a plain (without shortcuts) version of the algorithm.\n");
    fi;

    if test[1] = 1 and test[2] = 14 then
        Print(" The algorithm has not found an involution which behaves as \n");
        Print(" involution  of type t_n in PSp(2n,q). \n");
        Print(" This is NOT a definite answer, and it was given by  \n");
        Print(" a plain without shortcuts) version of the algorithm.\n");
        Print(" To be on the safe side, repeat the test.\n");
    fi;

    if l > 4 then
        Print(" The first run of the test has not produced the conclusive result,\n");
        Print(" the test was repeated ", test[l], " times.\n");
        Print(" The above conclusion is generated at the last run of the test.\n");
    fi;
    Print("\n");
    Print(" In the lists Test(G) and PlainTest(G), the first position is 0 or 1,\n");
    Print(" according to the answer\n");
    Print(" *symplectic* or *orthogonal* given by the test, the second position is\n");
    Print(" the rule number used in this decision, the third position is the actual type\n");
    Print(" of the involution used in the identification of the group; \n");
    Print(" (the latter was only used for the debugging of the algorithm).\n");
    Print(" The fourth number is the numbers of re-runs of the test\n");
    Print(" when the first run was inconclusive.\n");
end;


\end{verbatim}


\small

\begin{thebibliography}{99}



\bibitem{old} C.~Altseimer and A.~Borovik,
`Probabilistic  recognition of orthogonal and symplectic groups',
{\em Groups and Computations III}, (eds W.~M.Kantor and A.
Seress), de Gruyter, Berlin and New York, 2001, pp. 1--20.

\bibitem{artin}  E.~Artin, {\em Geometric Algebra}, Interscience Publishers,
New York, 1957.

\bibitem{car2} R.~W. ~Carter,  {\em Simple Groups of Lie Type}  (John Wiley
  and Sons, London a.o., 1972).

\bibitem{Carter-bk}  R. Carter, {\it Finite Groups of Lie Type:
Conjugacy Classes and Complex Characters}, Wiley, 1985.

\bibitem{lyons} D.~Gorenstein, R.~Lyons and R.~Solomon, `The
  classification of the finite simple groups 3',
 Mathematical Surveys and Monographs 40, AMS,
  Providence,  1998.

\bibitem{tay} D.~E.~Taylor, {\em The Geometry of the Classical
  Groups} (Heldermann, Berlin, 1992).

\bibitem{zs} K.~Zsigmondy, `Zur Theorie der
  Potenzreste',  {\em Montash. f\"ur Math.\ u.\ Phys.} (3) (1892) 265--284.


\end{thebibliography}
\end{document}